\documentclass[12pt]{article}

\usepackage{graphicx, xcolor}
\usepackage{epsfig}
\usepackage{amssymb}
\usepackage[T1]{fontenc}

\def\e{{\mathbf e}}

 \def\R{\mathbb{R}}

 \def\Z{{\bf Z}}

\def\C{\mathbb{C}}
\def\N{{\mathbb N}}

\newcommand{\proof}{\bf {Proof:} \rm}

\newtheorem{theorem}{Theorem}[section]

\newtheorem{lemma}{Lemma}[section]

\newtheorem{corollary}{Corollary}[section]
\newtheorem{definition}{Definition}[section]

\begin{document}

\title{Hermite Polynomials in Dunkl-Clifford Analysis}

\author{Minggang Fei\thanks{Corresponding author. E-mail: fei@ua.pt} \thanks{{\small School of Mathematical Sciences, University of Electronic Science and Technology of China, Chengdu, 610054, P. R. China}} \thanks{{\small Department of Mathematics, \textit{CIDMA - Center for Research and Development in Mathematics and Applications},
University of Aveiro, Aveiro, P-3810-193, Portugal}}, Paula
Cerejeiras$^{\ddagger}$ and Uwe K\"ahler$^{\ddagger}$}

\date{\today}

\maketitle

\begin{abstract} In this paper we present a generalization of the classical Hermite polynomials to the framework of Clifford-Dunkl operators. Several basic properties, such as orthogonality relations, recurrence formulae and associated differential equations, are established. Finally, an orthonormal basis for the Hilbert modules arising from the corresponding weight measures is studied. 
\end{abstract}



\section{Introduction}

It is well-known that classical harmonic analysis is linked to the invariance of the Laplacian under rotations. Unfortunately, many structures do not possess such invariance. In the 80's, C. Dunkl proposed a differential-difference operator associated to a given finite reflection group $W$. These operators are particularly adequate for the study of analytic structures with prescribed reflection symmetries, thus, providing a framework for a generalization of the classical theory of spherical harmonic functions (see \cite{Dunkl88}, \cite{Dunkl89}, \cite{DX01}, \cite{Roesler03}, \cite{CKR}, \cite{BCK}, \cite{FCK2009}, \cite{OSS}, etc.). These operators gained a renewed interest when it was realized that they had a physical interpretation, as they were naturally connected with certain Schr\"odinger operators for Calogero-Sutherland type quantum many body systems (see \cite{Roesler98}, \cite{Roesler03},\cite{FCK2010}, for more details). \\ 

In \cite{Roesler98}, R\"osler proposed a generalization of the classical Hermite polynomials systems to the multivariable case and proved some of their properties, such as Rodrigues and Mahler formulae and a generating relation,  analogies of the associated differential equations, together with its link to generalized Laguerre polynomials (see \cite{BF}). However, her generalization does not give a precise form for these polynomials. 

The study of special functions in the multivariable setting of Clifford analysis is not a new field. Already in his paper \cite{S88}, Sommen constructed a family of generalized Hermite polynomials by imposing axial symmetry and analysing the resulting Vekua-type system. By this technique he was successful  in obtaining the  orthogonality relation and a basis for the associated weighted $L_2$ space. His work proved to be the keystone for the multivariable generalizations of special functions within the Clifford analysis setting. In \cite{DeB}, De Bie used the approach developed in \cite{DeBS} for a further construction of such polynomials. Combining the previous technique of Sommen with a suitable Cauchy-Kovalevskaya extension he constructed concrete Clifford-Hermite polynomials of even degree. In fact, in the even case the powers of the Hermite operator are then scalar operators, thus making it easy to handle the  Dunkl-Laplace and -Euler operators. Unfortunately, no suggestion was made for handling the odd case.\\
 
It is the aim of this paper to complete De Bie's work by presenting the Clifford-Hermite polynomials of arbitrary positive degree related to the Dunkl operators. For that purpose, the authors will use the spherical representation formulae of the  Dunkl-Dirac operator obtained and studied in \cite{FCK2009}.\\

The paper is organized as follows. In Section 2 we collect the necessary basic facts regarding (universal) Clifford algebras and we present a spherical representation of Dunkl-Dirac operators. In Section 3 we present our main results. Namely, we give the definition of Clifford-Hermite polynomials related to the spherical representations of Dunkl operators for an arbitrary positive degree. Basic properties, such as orthogonality relations, recurrence formulae, and differential equations are proven. We finalize with the construction and study of the orthonormal basis for the Hilbert modules associated with the weight measures.

\section{Preliminaries}

\subsection{Clifford algebras}

Let $\e_1, \cdots, \e_d$ be an orthonormal basis of $\R^d$ satisfying the anti-commutation relationship 
$\e_i\e_j+\e_j\e_i=-2\delta_{ij}$, where $\delta_{ij}$ is the Kronecker symbol. One defines the universal real-valued Clifford algebra $\R_{0,d}$ as the $2^d$-dimensional associative algebra with basis given by $\e_0=1$ and $\e_A =\e_{h_1}\cdots \e_{h_n}$, where $A=\{(h_1, h_2, \cdots, h_n) : 1\leq h_1<h_2<\cdots<h_n\leq d \} $. Hence, each element $x\in\R_{0,d}$ can be written as $x=\sum_Ax_A\e_A$, $x_A\in\R$. In what follows, $sc[x]=x_0$ will denote the scalar part of $x\in\R_{0,d}$, while a vector $(x_1, x_2,\cdots,x_d) \in \R^d$ will be identified with the element $x=\sum_{i=1}^dx_i\e_i$. \\

We define the Clifford conjugation as a linear action from $\R_{0,d}$ into itself, which acts on the basis elements as  $$\bar{1} = 1,~~~\bar{\e}_i=-\e_i, ~i= 1, \cdots, d$$ and possess the anti-involution property $\overline{\e_i\e_j}=\bar{\e}_j\bar{\e}_i.$ An important property of $\R_{0,d}$ is that each non-zero vector $x\in\R^d$ has a multiplicative inverse given by $x^{-1}=\frac{\bar{x}}{\|x\|^2}=\frac{-x}{\|x\|^2},$ where the norm $\| \cdot \|$ is the usual Euclidean norm. 

Therefore, in Clifford notation, the reflection $\sigma_\alpha x$ of a vector $x\in\R^d$ with respect to the hyperplane $H_{\alpha}$ orthogonal to a given $\alpha\in\R^d\backslash\{0\},$ is $$\sigma_\alpha
x=-\alpha x \alpha^{-1} = x + \frac{2 \langle x, \alpha \rangle }{\| \alpha \|^2} \alpha,$$ with $\langle \cdot, \cdot \rangle$ denoting the standard Euclidean inner product. \\

Functions spaces are introduced as follows. A $\C \otimes \R_{0,d}$-valued function $f$ in an open set $\Omega \subset \R^d $ has a representation $f=\sum_A\e_A f_A$, with components $f_A: \Omega\rightarrow\C$. Function spaces of Clifford-valued functions are established as modules over $\R_{0,d}$ by imposing its coefficients $f_A$ to be in the corresponding real-valued function space. For example, $f=\sum_A\e_A f_A \in L_2(\Omega; \C \otimes \R_{0,d})$ if and only $f_A \in L_2(\Omega), \forall A.$ When no ambiguity arises, we will use the complex valued notation for the correspondent Clifford-valued module.\\ 

\subsection{Dunkl operators in Clifford setting}

A finite set $R\subset\R^d\backslash\{0\}$ is called a root system if $R\bigcap \alpha \R^d =\{\alpha,-\alpha\}$ and $\sigma_{\alpha}R=R$ for all $\alpha\in R$. For a given root system $R$ the set of reflections $\sigma_{\alpha}$, $\alpha\in R,$ generates a finite group $W\subset O(d)$, called the finite reflection group (or Coxeter group) associated with $R$. All reflections in $W$ correspond to suitable pairs of roots. For a given $\beta\in\R^d\backslash\bigcup_{\alpha\in R}H_{\alpha}$, we fix the positive subsystem $R_{+}=\{\alpha\in R|\langle\alpha,\beta\rangle>0\}$, i.e. for each $\alpha\in R$ either $\alpha\in R_{+}$ or $-\alpha\in R_{+}$.\\

A function $\kappa: R\rightarrow\C$ is called a multiplicity function on the root system if it is invariant under the action of the associated reflection group $W$. This means that $\kappa$ is constant on the conjugacy classes of reflections in $W$. For abbreviation, we introduce the index $\gamma_{\kappa}=\sum_{\alpha\in R_{+}}\kappa(\alpha)$ and the Dunkl-dimension $\mu=2\gamma_{\kappa}+d.$

For each fixed positive subsystem $R_+$ and multiplicity function $\kappa$ we have, as invariant operators, the differential-difference operators (also called Dunkl operators):
\begin{eqnarray}
T_if(x)=\frac{\partial}{\partial x_i}f(x)+\sum_{\alpha\in R_{+}} \kappa(\alpha)\frac{f(x)-f(\sigma_{\alpha}x)}{\langle\alpha,x\rangle}\alpha_i, \qquad i=1,  \cdots, d,
\end{eqnarray}
for $f\in C^1(\R^d)$. In the case of $\kappa=0$, the operators coincide with the corresponding partial derivatives.  Therefore, these differential-difference operators can be regarded as the equivalent of partial derivatives related to given finite reflection groups. More important, these operators commute, that is, $T_iT_j=T_jT_i$.\\

In this paper we will assume $Re (\kappa) \geq 0$ and $\gamma_{\kappa}>0$. Based on these real-valued operators we introduce the Dunkl-Dirac operator in $\R^d$ associated to the reflection group $W,$ and multiplicity function $\kappa,$ as (\cite{CKR},\cite{OSS})
\begin{eqnarray}
D_hf=\sum_{i=1}^d\e_iT_if.
\end{eqnarray} As in the classic case, the Dunkl-Dirac operator factorize the Dunkl Laplacian 
in $\R^d$  by $$\Delta_h=-D_h^2=\sum_{i=1}^dT_i^2.$$ Functions belonging to the kernel of Dunkl-Dirac operator
will be called Dunkl-monogenic functions. As usual, functions belonging to be the kernel of Dunkl Laplacian will be called Dunkl-harmonic functions.

For the construction of Hermite polynomials of arbitrary positive degree we require the following two lemmas regarding the decomposition into spherical coordinates $x = r\omega, r=|x|,$ of the Dunkl-Dirac operator. 

\begin{lemma}[Theorem 3.1 in \cite{FCK2009}]\label{le:2.1} In spherical coordinates the Dunkl-Dirac operator has the following form:
\begin{eqnarray}
D_hf(x) =\omega \left(   \partial_r+\frac{1}{r} \Gamma_{\kappa} \right)f(x) = \omega \left[   \partial_r+\frac{1}{r}  \left(   \gamma_{\kappa}+\Phi_{\omega}+\Psi \right) \right] f(r\omega),
\end{eqnarray}
where
\begin{eqnarray*}
\Phi_{\omega}f(x) =-\sum_{i<j}\e_i\e_j(x_i\partial_{x_j}-x_j\partial_{x_i})f(x),
\end{eqnarray*}
and
\begin{eqnarray*}
\Psi f(x)=-\sum_{i<j}\e_i\e_j\sum_{\alpha\in
R^+}\kappa(\alpha)\frac{f(x)-f(\sigma_{\alpha}x)}{\langle\alpha,x\rangle}(x_i\alpha_j-x_j\alpha_i)-\sum_{\alpha\in
R^+}\kappa(\alpha)f(\sigma_{\alpha}x),
\end{eqnarray*}
for $f\in C^1(\R^d)$.
\end{lemma}

\begin{lemma}[Theorems 3.2 and 3.3 in \cite{FCK2009}] \label{le:2.2} 
The operator $\Gamma_{\omega}$ satisfies
\begin{enumerate}
\item $ \Gamma_{\omega} f(r)=0,$ if $f$ is a radial function.
\item $\Gamma_{\omega}(\omega)=(\mu-1)\omega,$
\item $\Gamma_{\omega}P_n(\omega)=-nP_n(\omega),$ 
\item $\Gamma_{\omega}(\omega P_n(\omega))=(\mu+n-1)\omega
P_n(\omega).$
\end{enumerate} where $P_n$ denotes a homogeneous Dunkl-monogenic function of degree $n\in\Z.$ 
\end{lemma}

Henceforward, we denote by $M_n$ the space of all homogeneous Dunkl-monogenic polynomials of degree $n \in\N. $
We have then 
\begin{lemma} \label{le:2.3} Let $s\in\N$ and $P_n\in M_n$. Then for any radial function $f(r)=f(|x|)$ it is valid
\begin{enumerate}
\item $ D_h(f(r)P_n(x))=\omega f'(r)P_n(x),$
\item $D_h(\omega f(r)P_n(x))=-\left(f'(r)+\frac{\mu+2n-1}{r}P_n(x)\right)$
\item $D_h(x^sP_n(x))=\left\{
       \begin{array}{ll}
        -sx^{s-1}P_n(x), \ \  s \ even,\\
        \ \\
        -(s+\mu+2n-1)x^{s-1}P_n(x), \ \  s \ odd.
        \end{array}
        \right.$
\end{enumerate}
\end{lemma}

\section{Hermite Polynomials in Dunkl-Clifford Analysis}

We denote by $L^2(\R^d; e^{x^2})$ the weighted $L^2$-space of Clifford-valued measurable functions in $\R^d$ induced by the inner product
$$(f,g)_H=\int_{\R^d}\overline{f(x)}g(x)e^{x^2}h_{\kappa}^2(x)dx.$$
We remark that $L^2(\R^d;e^{x^2})$ is a right Hilbert module over $\C \otimes \R_{0,d}$.

For our purpose, it is required to analyse the behaviour of the inner product for functions of type  $f(x)=x^sP_n(x)$, where $P_n\in M_n.$ 

\begin{lemma}\label{le:3.1} If we let $n,s,t\in\N$ and $P_n\in M_n$, then
\begin{eqnarray*}
(x^sP_n, x^tP_n)_H=\left\{
       \begin{array}{cl}
        (-1)^{\frac{s+t}{2}}\frac{1}{2}\Gamma(\frac{s+t+2n+\mu}{2})\|P_n\|_{\kappa}^2 & {\rm , ~if } ~ s ~{\rm and} ~t ~{\rm are~even,}\\
        \ \\
        (-1)^{\frac{s+t}{2}+1}\frac{1}{2}\Gamma(\frac{s+t+2n+\mu}{2})\|P_n\|_{\kappa}^2 &  {\rm , ~if } ~ s ~{\rm and} ~t ~{\rm are~odd,}\\
        \ \\
        0 & {\rm , ~if} ~ s ~{\rm and} ~t ~{\rm have~different~parity},  
        \end{array}
        \right.
\end{eqnarray*}
where $\|P_n\|_{\kappa}=(\int_{S^{d-1}}|P_n(\omega)|^2h_{\kappa}^2(\omega)d\Sigma(\omega))^{1/2}$ is the usual spherical norm of $P_n$ in Dunkl analysis.
\end{lemma}

\proof Using the spherical coordinates $x=r\omega$, $r=|x|$, we have, 
\begin{eqnarray*}
(x^sP_n, x^tP_n)_H & = &  \int_{\R^d}\overline{P_n(x)}\bar{x}^sx^tP_n(x)e^{x^2}h_{\kappa}^2(x)dx\\
         & = & \int_0^{\infty}r^nr^sr^tr^ne^{r^2}r^{2\gamma_{\kappa}}r^{d-1}dr \int_{S^{d-1}}\overline{P_n(\omega)}\bar{\omega}^s\omega^tP_n(\omega)h_{\kappa}^2(\omega)d\Sigma(\omega)\\
        & = & \frac{1}{2}\Gamma(\frac{s+t+2n+\mu}{2})\int_{S^{d-1}}\overline{P_n(\omega)}\bar{\omega}^s\omega^tP_n(\omega)h_{\kappa}^2(\omega)d\Sigma(\omega).
\end{eqnarray*}

First, we consider the case in which both $s$ and $t$ are even. 
Let $s=2a$ and $t=2b,$ for some  $a,b\in\N.$ Then
\begin{eqnarray*}
(x^sP_n, x^tP_n)_H & = &  \frac{1}{2}\Gamma(\frac{s+t+2n+\mu}{2})(-1)^{a+b}\int_{S^{d-1}}\overline{P_n(\omega)}P_n(\omega)h_{\kappa}^2(\omega)d\Sigma(\omega)\\
                 & = & (-1)^{\frac{s+t}{2}}\frac{1}{2}\Gamma(\frac{s+t+2n+\mu}{2})\|P_n\|_{\kappa}^2.
\end{eqnarray*}

In a similar way, we obtain 
$$  (x^sP_n, x^tP_n)_H =  (-1)^{\frac{s+t}{2} +1}\frac{1}{2}\Gamma(\frac{s+t+2n+\mu}{2})\|P_n\|_{\kappa}^2 $$ 
when both $s$ and $t$ are odd. 

Now, when $s=2a$ is even and $t=2b+1$ is odd, with $a,b\in\N$, we get 
$$
(x^sP_n, x^tP_n)_H = \frac{1}{2}\Gamma(\frac{s+t+2n+\mu}{2})(-1)^{a+b}\int_{S^{d-1}}\overline{P_n(\omega)}\omega 
P_n(\omega)h_{\kappa}^2(\omega)d\Sigma(\omega).
$$ If $P_n\in M_n$ we have that $xP_n(x)$ is a homogeneous Dunkl-harmonic polynomial of degree $n+1$ (see \cite{DX01}, Lemma 5.1.10). Hence, by the orthogonality property of Dunkl-harmonics of different degree, we obtain
\begin{eqnarray*}
\int_{S^{d-1}}\overline{P_n(\omega)}\omega
P_n(\omega)h_{\kappa}^2(\omega)d\Sigma(\omega)=0, 
\end{eqnarray*}
so that $(x^sP_n, x^tP_n)_H=0.$ The remaining case is analogous. $\qquad \blacksquare$\\




Following \cite{DeBS}, we now introduce the vector space
$$R(P_n) = \left\{ \sum_{j=0}^m a_jx^jP_n(x) | m\in\N, a_j\in \C, P_n\in M_n \right\}.$$
In particular, we have $R(1) = \left\{ \sum_{j=0}^m a_jx^j | m\in\N, a_j\in\C \right\}.$

Also, we introduce the operator $D_{+}=D_h-2x.$ It is easy to see that $D_h(R(P_n))\subset R(P_n),$ due to Lemma \ref{le:2.3}. Hence, the following properties of the inner product $(\cdot,\cdot)_H$ are valid.

\begin{lemma}\label{le:3.2} For fixed $P_n \in M_n$ it holds
$$(D_{+}(p P_n), q P_n)_H=( p P_n,D_h( q P_n))_H,$$
 where $p, q \in R(1).$
\end{lemma}

\proof It suffices to prove that
\begin{enumerate}
\item $ (D_{+}(x^{2s}P_n),x^{2t}P_n)_H=(x^{2s}P_n,D_h(x^{2t}P_n))_H;$
\item $(D_{+}(x^{2s+1}P_n),x^{2t+1}P_n)_H=(x^{2s+1}P_n,D_h(x^{2t+1}P_n))_H;$
\item $(D_{+}(x^{2s}P_n),x^{2t+1}P_n)_H=(x^{2s}P_n,D_h(x^{2t+1}P_n))_H;$
\item $(D_{+}(x^{2s+1}P_n),x^{2t}P_n)_H=(x^{2s+1}P_n,D_h(x^{2t}P_n))_H.$
\end{enumerate}
The first two identities are immediate since
\begin{eqnarray*}
(D_{+}(x^{2s}P_n),x^{2t}P_n)_H = (x^{2s}P_n,D_h(x^{2t}P_n))_H & = & 0,\\
(D_{+}(x^{2s+1}P_n),x^{2t+1}P_n)_H = (x^{2s+1}P_n,D_h(x^{2t+1}P_n))_H& = & 0,
\end{eqnarray*}
by our Lemma \ref{le:3.1}. Identities $3.$ and $4.$ can be proved in a similar way.

Now, on one hand, we have 
\begin{eqnarray*}
(D_{+}(x^{2s}P_n),x^{2t+1}P_n)_H & = & -2s(D_{+}(x^{2s-1}P_n),x^{2t+1}P_n)_H-2(x^{2s+1}P_n,x^{2t+1}P_n)_H\\
& = & -2s(-1)^{\frac{2s+2t}{2}+1}\frac{1}{2}\Gamma(\frac{2s+2t+2n+\mu}{2})\|P_n\|_{\kappa}^2\\
& &\qquad -2(-1)^{\frac{2s+2t+2}{2}+1}\frac{1}{2}\Gamma(\frac{2s+2t+2n+\mu}{2}+1)\|P_n\|_{\kappa}^2\\
& = &2s(-1)^{s+t}\frac{1}{2}\Gamma(\frac{2s+2t+2n+\mu}{2})\|P_n\|_{\kappa}^2\\
& &\qquad -2(-1)^{s+t}\frac{1}{2}\Gamma(\frac{2s+2t+2n+\mu}{2}+1)\|P_n\|_{\kappa}^2.
\end{eqnarray*}
On the other hand,
\begin{eqnarray*}
(x^{2s}P_n,D_h(x^{2t+1}P_n))_H & = & -(2t+1+2n+\mu-1)(x^{2s}P_n,x^{2t}P_n)_H\\
& = & -(2t+2n+\mu)(-1)^{\frac{2s+2t}{2}}\frac{1}{2}\Gamma(\frac{2s+2t+2n+\mu}{2})\|P_n\|_{\kappa}^2\\
& = & 2s(-1)^{s+t}\frac{1}{2}\Gamma(\frac{2s+2t+2n+\mu}{2})\|P_n\|_{\kappa}^2\\
& &\qquad -2(-1)^{s+t}\frac{1}{2}(\frac{2s+2t+2n+\mu}{2})\Gamma(\frac{2s+2t+2n+\mu}{2})\|P_n\|_{\kappa}^2\\
& = & 2s(-1)^{s+t}\frac{1}{2}\Gamma(\frac{2s+2t+2n+\mu}{2})\|P_n\|_{\kappa}^2\\
& &\qquad 2(-1)^{s+t}\frac{1}{2}\Gamma(\frac{2s+2t+2n+\mu}{2}+1)\|P_n\|_{\kappa}^2.
\end{eqnarray*}
From these two relations one gets
\begin{eqnarray*}
(D_{+}(x^{2s}P_n),x^{2t+1}P_n)_H & = & (x^{2s}P_n,D_h(x^{2t+1}P_n))_H.
\end{eqnarray*}
This completes the proof. $\qquad \blacksquare$

We now recall the definition of Hermite polynomials in Dunkl-Clifford analysis.

\begin{definition}\label{def:3.1} Fix $P_n\in M_n.$ Then, for each $s \in \N_0$
$$H_{s, \mu, P_n}(x):=(D_{+})^sP_n(x)$$
is a Dunkl-Clifford-Hermite polynomial of degree $(s,n).$
\end{definition}

\textbf{Remark} \textit{Dunkl-Clifford-Hermite polynomials depend on the initial choice of the monogenic polynomial $P_n.$}\\

Due to Lemma \ref{le:2.1}, we can now apply the definition to the case of the Hermite polynomials of an arbitrary positive degree. In fact, due to this lemma, we have  
$$H_{s,\mu, P_n}(x)=H_{s,\mu, 1}(x)P_n(x),$$
where $H_{s,\mu, 1} \in R(1)$ depends only on the degree $s.$ So, it is easy to conclude that $H_{s,\mu, P_n} \in R(P_n)$.

We give here the explicit form of the first Dunkl-Clifford-Hermite polynomials.

\begin{eqnarray*}
H_{0,\mu, P_n}(x) & = & P_n(x),\\
H_{1,\mu, P_n}(x) & = & -2xP_n(x),\\
H_{2,\mu, P_n}(x) & = & [4x^2+2(\mu+2n)]P_n(x),\\
H_{3,\mu, P_n}(x) & = & -[8x^3+4(\mu+2n+2)x]P_n(x),\\
H_{4,\mu, P_n}(x) & = & [16x^4+16(\mu+2n+2)x^2+4(\mu+2n+2)(\mu+2n)]P_n(x),\\
& \vdots &
\end{eqnarray*}

Using this definition, we obtain a straightforward recurrence relation.

\begin{lemma}(Recurrence relation) \label{le:3.3} For each fixed $P_n \in M_n,$ the recurrence relation 
$$H_{s,\mu, P_n}(x) = D_{+}H_{s-1,\mu, P_n}(x), ~~s\in \N,$$ holds.
\end{lemma}

Also, we can prove a Rodrigues' formula in the general case for Dunkl-Clifford-Hermite polynomials of arbitrary positive degree.

\begin{theorem}(Rodrigues' formula) \label{th:3.1} $H_{s,\mu}(P_n)(x)$ is also determined by
$$H_{s,\mu, P_n}(x)=e^{r^2}(D_h)^s(e^{-r^2}P_n(x)), ~~|x| = r.$$
\end{theorem}

\proof The key point in our proof is the following identity relating the Dunkl-Dirac operator $D_h$ with the $D_+$ operator. For any $f\in C^1(\R^d)$, we have
\begin{eqnarray}
e^{r^2}D_h(e^{-r^2}f)&=&e^{r^2}[\omega(\partial_r+\frac{1}{r}\Gamma_{\omega})](e^{-r^2}f)  \nonumber \\
                    &=&e^{r^2}[\omega(e^{-r^2}(-2r)f+e^{-r^2}\partial_rf+\frac{1}{r}e^{-r^2}\Gamma_{\omega}f)]  \nonumber \\
                    &=&-2xf+D_hf  \nonumber \\
                    &=&D_{+}f.  \label{result1}
\end{eqnarray}
Therefore, 
\begin{eqnarray*}
e^{r^2}(D_h)^s(e^{-r^2}P_n(x))&=&e^{r^2}(D_h)^{s-1}(e^{-r^2}e^{r^2}D_h(e^{-r^2}P_n(x)))\\   \nonumber \\
                              &=&e^{r^2}(D_h)^{s-1}(e^{-r^2}D_{+}P_n(x))
\end{eqnarray*} Proceeding  recursively we obtain
$$ e^{r^2}(D_h)^s(e^{-r^2}P_n(x)) = (D_{+})^{s}P_n(x) = H_{s,\mu, P_n}(x). \qquad \blacksquare$$

The orthogonality between Dunkl-Clifford-Hermite polynomials 
is expressed as follows.

\begin{lemma}(Orthogonality relation) \label{le:3.4} If $s\neq t$, then
$$(H_{s,\mu, P_n},H_{t,\mu, P_n})_H=0.$$
\end{lemma}

Again, the proof of the orthogonality is rather straightforward. It relays on the fact that $H_{s, \mu, P_n} = D^s_+ P_n \in R(P_n),$ on applying  Lemma \ref{le:3.2} for interchanging $D_+^s$ with $D_h^s,$ and using Lemma \ref{le:2.3}, property 3. to conclude that $D_h^s (H_{t, \mu, P_n}) = 0$ whenever $t <s.$  

\begin{corollary}\label{co:3.1} For every fixed $P_n \in M_n$ the polynomials $H_{s,\mu, P_n}, ~s\in \N_0$, forms a  basis of $R(P_n)$.
\end{corollary}

We are now in a position to prove that Dunkl-Clifford-Hermite polynomials satisfy a differential equation in Dunkl case. This equation is given as follows.

\begin{theorem}(Differential equation) \label{th:3.2} For each fixed $P_n \in M_n,$ the Dunkl-Clifford-Hermite polynomial $H_{s,\mu, P_n}$ satisfies the differential equation
$$D_h^2H_{s,\mu, P_n}-2xD_hH_{s,\mu, P_n}-C(s,\mu,n)H_{s,\mu, P_n} =0,$$
where
$$
C(s,\mu,n) = \left\{
       \begin{array}{ll}
        2s, \ \ if ~ s ~ even,\\
        \ \\
        2(s+\mu+2n-1), \ \ if ~ s ~ odd.\\
        \end{array}
        \right.$$
\end{theorem}

\proof The proof relays on the fact that $H_{s,\mu, P_n} = H_{s,\mu, 1} P_n,$ with $H_{s,\mu, 1} \in R(1).$ Hence, when one applies the Dunkl operator to $H_{s,\mu, P_n}$ it reduce the degree of the polynomial $H_{s,\mu, 1}$  by $1$ (by Lemma \ref{le:2.3}), that is, it exists a polynomial $p$ of degree $s-1$ such that $D_h H_{s,\mu, P_n} = p P_n.$ 

Now, since the polynomials $H_{s,\mu, 1}, ~s\in \N_0$, forms a  basis of $R(1)$ (Corollary \ref{co:3.1}) we can write 
$$D_h H_{s,\mu, P_n} = p P_n = \left( \sum_{j=0}^{s-1} b_j H_{j, \mu, 1} \right) P_n = \sum_{j=0}^{s-1} b_j H_{j, \mu, P_n} ,$$ for some $b_0, b_1, \cdots, b_{s-1} \in \C.$

For $0 \leq i < s-1,$ we consider the inner product $(H_{i,\mu, P_n},\sum_{j=0}^{s-1} b_j H_{j,\mu, P_n})_H. $ On one hand, 
$$(H_{i,\mu, P_n},\sum_{j=0}^{s-1} b_j H_{j,\mu, P_n})_H = b_i \|H_{i,\mu, P_n} \|_H^2. $$
On the other hand, 
\begin{eqnarray*}
(H_{i,\mu, P_n},\sum_{j=0}^{s-1} b_j H_{j,\mu, P_n})_H & = &  (H_{i,\mu, P_n}, D_h H_{s,\mu, P_n})_H \\
& = & (D_+ H_{i,\mu, P_n}, H_{s,\mu, P_n})_H\\
&=&(H_{i+1,\mu, P_n} , H_{s,\mu, P_n})_H\\
&=& 0.
\end{eqnarray*} These both conditions imply each $b_i =0, ~i=0, 1, \cdots, s-2,$ so that 
\begin{equation}
D_h H_{s,\mu, P_n}  =  \sum_{j=0}^{s-1} b_j H_{j, \mu, P_n} = b_{s-1} H_{s-1, \mu, P_n}. \label{Eq:1}
\end{equation}

We set $C(s,\mu,n)=b_{s-1}.$ 

On one hand, by applying the $D_+$ operator on both sides of (\ref{Eq:1}), we obtain  
\begin{equation}
D_+D_hH_{s,\mu, P_n} = C(s,\mu,n)D_+H_{s-1,\mu, P_n} = C(s,\mu,n) H_{s,\mu, P_n}. \label{eq:1}
\end{equation}

On the other hand, due to (\ref{result1}) we have 
\begin{eqnarray}
D_+D_hH_{s,\mu, P_n } &=&e^{r^2}D_h(e^{-r^2}D_h H_{s,\mu, P_n} )  \nonumber \\
                       &=&e^{r^2}\omega(\partial_r+\frac{1}{r}\Gamma_{\omega})(e^{-r^2}D_hH_{s,\mu, P_n} ) \nonumber \\
                       &=&-2r\omega D_hH_{s,\mu, P_n} +D_h^2H_{s,\mu, P_n} \nonumber \\
                       &=&-2xD_hH_{s,\mu, P_n}+D_h^2H_{s,\mu, P_n}. \label{eq:2}
\end{eqnarray}
Combining (\ref{eq:1}) and (\ref{eq:2}) we get
\begin{eqnarray}
D_h^2H_{s,\mu, P_n}-2x D_hH_{s,\mu, P_n} =C(s,\mu,n)H_{s,\mu, P_n}. \label{eq:3}
\end{eqnarray}

Finally, taking into account that $H_{s,\mu, P_n} = \sum_{j=0}^s a_j x^j P_n,$ and Lemma \ref{le:2.3} then equality (\ref{eq:3}) yields

\begin{eqnarray*}
& &D_h^2H_{s,\mu, P_n}(x)-2xD_hH_{s,\mu, P_n}(x)\\
& & \\
& &\qquad=\left\{
       \begin{array}{ll}
        2sa_sx^sP_n(x)+{\rm ~terms ~of ~lower~ order}, \ \  {\rm ~if ~} s {\rm ~ even},\\
        \ \\
        2(s+\mu+2n-1)a_sx^sP_n(x)+{\rm ~terms ~of ~lower~ order}, \ \ {\rm ~if ~} s {\rm ~ odd}.\\
        \end{array}
        \right.
\end{eqnarray*}
Comparing the coefficients of the highest terms on both sides of (\ref{eq:3}) gives
\begin{eqnarray*}
C(s,\mu,n)=\left\{
       \begin{array}{ll}
        2s, \ \  {\rm ~if ~} s {\rm ~ even},\\
        \ \\
        2(s+\mu+2n-1), \ \  {\rm ~if ~} s {\rm ~ odd}.\\
        \end{array}
        \right.
\end{eqnarray*}
This completes the proof.  $\qquad \blacksquare$\\


\begin{lemma}(Three terms recurrence) \label{le:3.5} For a fixed $P_n \in M_n$ and $s \in \N$ we have
$$H_{s+1,\mu, P_n} = -2x H_{s,\mu, P_n} +C(s,\mu,n) H_{s-1,\mu, P_n}.$$
\end{lemma}

\proof In fact,  \begin{eqnarray*} 
H_{s+1,\mu, P_n}& = & D_+H_{s,\mu, P_n} \\
                &=&(D_h-2x)H_{s,\mu, P_n}\\
                &=&-2x H_{s,\mu, P_n}+ C(s,\mu,n) H_{s-1,\mu, P_n}. \qquad \blacksquare
\end{eqnarray*}

\begin{corollary} \label{co:3.2} From the three terms recurrence formula we get
$$H_{s,\mu, P_n} = \left\{ \begin{array}{ll} 
\sum_{j=0}^t a^{2t}_{2j}x^{2j}P_n, & {\rm ~if ~} s=2t\\
 & \\
 \sum_{j=0}^t a^{2t+1}_{2j+1}x^{2j+1}P_n, & {\rm ~if ~} s=2t+1
\end{array} \right. .$$
\end{corollary}

Furthermore, as we have $H_{s,\mu, P_n} = H_{s,\mu,1} P_n,$ with $H_{s,\mu,1} \in R(1),$ we can use the recurrence relation (Lemma \ref{le:3.1}) together with the differential equation (Theorem \ref{th:3.2}) in order to compare the Dunkl-Clifford-Hermite polynomials $H_{s,\mu,P_n}$ with orthogonal polynomials on the real line. 

\begin{theorem}\label{th:3.3} For each fixed $P_n \in M_n$ and $s \in \N_0$ we have
\begin{eqnarray*}
H_{s,\mu,n}(x)=\left\{
       \begin{array}{ll}
        2^s(\frac{s}{2})! \ L_{\frac{s}{2}}^{\frac{\mu}{2}+n-1}(|x|^2), & {\rm ~if ~} s {\rm ~ even}\\
        & \\
        -2^s(\frac{s-1}{2})! \ x \  L_{\frac{s-1}{2}}^{\frac{\mu}{2}+n}(|x|^2), & {\rm ~if ~} s {\rm ~ odd}.\\
        \end{array}
        \right.
\end{eqnarray*}
where $L_s^{\alpha}(x)=\sum_{j=0}^s\frac{\Gamma(s+\alpha+1)}{j!(s-j)!\Gamma(j+\alpha+1)}(-x)^j$ denotes the generalized Laguerre polynomial on the real line.
\end{theorem}

\proof From Corollary \ref{co:3.2}, Lemmas \ref{le:3.5} and \ref{le:2.3}, we obtain the following relation between the coefficients of an arbitrary Dunkl-Clifford-Hermite polynomial 
\begin{eqnarray}
\left\{
       \begin{array}{ll}
        a_{2j}^{2t}=2(j+1)(2j+\mu+2n)a_{2j+2}^{2t-2}+2(4j+\mu+2n)a_{2j}^{2t-2}+4a_{2j-2}^{2t-2},\\
        \ \\
        a_{2j+1}^{2t+1}=2(j+1)(2j+\mu+2n+2)a_{2j+3}^{2t-1}+2(4j+\mu+2n+2)a_{2j+1}^{2t-1}+4a_{2j-1}^{2t-1}.\\
        \end{array}
        \right. \label{Eqn:1}
\end{eqnarray} 

Using Theorem \ref{th:3.2} and Lemma \ref{le:2.3} we obtain
\begin{eqnarray}
\left\{
       \begin{array}{ll}
        2j(2j+\mu+2n-2)a_{2j}^{2t}=4(t-j+1)a_{2j-2}^{2t},\\
        \ \\
        2j(2j+\mu+2n)a_{2j+1}^{2t+1}=4(t-j+1)a_{2j-1}^{2t+1}\\
        \end{array}
        \right. .\label{Eqn:2}
\end{eqnarray}
From  (\ref{Eqn:2}) we obtain
\begin{eqnarray}
\left\{
       \begin{array}{ll}
       a_{2j}^{2t}=\frac{t-j+1}{j(j+\frac{\mu}{2}+n-1)}a_{2j-2}^{2t}=\cdots=\frac{t!}{j!(t-j)!}\frac{\Gamma(\frac{\mu}{2}+n)}{\Gamma(\frac{\mu}{2}+n+j)}a_0^{2t},\\
        \ \\
        a_{2j+1}^{2t+1}=\frac{t-j+1}{j(j+\frac{\mu}{2}+n)}a_{2j-1}^{2t+1}=\cdots=\frac{t!}{j!(t-j)!}\frac{\Gamma(\frac{\mu}{2}+n+1)}{\Gamma(\frac{\mu}{2}+n+j+1)}a_1^{2t+1}.\\
        \end{array}
        \right.
\end{eqnarray}

Using equalities (\ref{Eqn:1}) and (\ref{Eqn:2}) again we have
\begin{eqnarray}
\left\{
       \begin{array}{ll}
       a_{0}^{2t}=2^2(\frac{\mu}{2}+n+t-1)a_{0}^{2t-2}=\cdots=2^{2t}\frac{\Gamma(\frac{\mu}{2}+n+t)}{\Gamma(\frac{\mu}{2}+n)}a_0^{0}=2^{2t}\frac{\Gamma(\frac{\mu}{2}+n+t)}{\Gamma(\frac{\mu}{2}+n)},\\
        \ \\
        a_{1}^{2t+1}=2^2(\frac{\mu}{2}+n+t)a_{1}^{2t-1}=\cdots=2^{2t}\frac{\Gamma(\frac{\mu}{2}+n+t+1)}{\Gamma(\frac{\mu}{2}+n+1)}a_1^{1}=2^{2t}\frac{\Gamma(\frac{\mu}{2}+n+t+1)}{\Gamma(\frac{\mu}{2}+n+1)}(-2).\\
        \end{array}
        \right.
\end{eqnarray}
Comparing with the definition of the generalized Laguerre polynomials yields the results of the theorem. $\qquad \blacksquare$

Finally, if we let $\{P_n^{(j)}|j=1,\cdots,\left(\begin{array} {cc} n+d-2 \\ n
\end{array} \right)\}$ be an orthonormal basis of $M_n$, i.e., $\frac{1}{|S^{d-1}|}\int_{S^{d-1}}\overline{P^{(i)}_n(\omega)}P^{(j)}_n(\omega)h_{\kappa}^2(\omega)d\Sigma(\omega)=\delta_{ij}$, then using the method introduced in \cite{DSS} it holds

\begin{theorem}\label{th:3.4} The set $\left\{ \frac{H_{s,\mu, P^{(j)}_n}}{\sqrt{\gamma_{s,\mu,n}}}|s,n,j\in\N, j\leq\left(\begin{array} {cc} n+d-2 \\ n
\end{array} \right)  \right\}$ is an orthonormal basis for $L^2(\R^d;e^{x^2})$, where $\gamma_{s,\mu,n}$ is given by
\begin{eqnarray*}
\gamma_{s,\mu,n}&=&(H_{s,\mu, P^{(j)}_n} ,H_{s,\mu, P^{(j)}_n} )_H\\
                &=&\left\{
       \begin{array}{ll}
        4^s(\frac{s}{2})!\pi^{\frac{d}{2}}\frac{\Gamma(\frac{s+\mu}{2}+n)}{\Gamma(\frac{d}{2})}, \ \  s \ even,\\
        \ \\
        4^s(\frac{s-1}{2})!\pi^{\frac{d}{2}}\frac{\Gamma(\frac{s+\mu+1}{2}+n)}{\Gamma(\frac{d}{2})}, \ \  s \ odd.\\
        \end{array}
        \right.
\end{eqnarray*}
\end{theorem}

\proof  We use the method described in \cite{DSS} to show that $\{ H_{s,\mu, P^{(j)}_n} \}$ is an orthogonal basis of
$L^2(\R^d;e^{x^2})$, here we only calculate the normalization
constants $\gamma_{s,\mu,n}$, that is
\begin{eqnarray*}
\gamma_{s,\mu,n}&=&(H_{s,\mu}(P^{(j)}_n)(x),H_{s,\mu}(P^{(j)}_n)(x))_H\\
&=&\frac{1}{C(s,\mu,n)}(D_+D_hH_{s,\mu}(P^{(j)}_n)(x),H_{s,\mu}(P^{(j)}_n)(x))_H\\
                                                 &=&\frac{1}{C(s,\mu,n)}(D_hH_{s,\mu}(P^{(j)}_n)(x),D_hH_{s,\mu}(P^{(j)}_n)(x))_H\\
                                                 &=&C(s,\mu,n)(H_{s-1,\mu}(P^{(j)}_n)(x),H_{s-1,\mu}(P^{(j)}_n)(x))_H\\
                                                 &=&C(s,\mu,n)C(s-1,\mu,n)\cdots C(1,\mu,n)(P^{(j)}_n)(x),P^{(j)}_n)(x))_H\\
                                                 &=&C(s,\mu,n)C(s-1,\mu,n)\cdots C(1,\mu,n)\frac{1}{2}\Gamma(\frac{\mu}{2}+n)\frac{2\pi^{\frac{d}{2}}}{\Gamma(\frac{d}{2})}.
\end{eqnarray*}
Substituting the coefficients $C(s,\mu,n)$ by their exact values
gives the desired formulae. $\qquad \blacksquare$

\bigskip
\bigskip
\bigskip

\noindent{\large \bf  Acknowledgements}

\bigskip

\noindent The authors were (partially) supported by {\it CIDMA - Centro de Investiga\c c\~ao e Desenvolvimento em Matem\'atica e Aplica\c c\~oes} of the University of Aveiro. The first author is the recipient of a grant from
 {\it Funda\c{c}\~{a}o para Ci\^{e}ncia e a Tecnologia (Portugal)} with grant No.: SFRH/BPD/41730/2007.

 \end{document}